\documentclass[12pt]{article}
\usepackage{amsmath}
\usepackage{amsfonts}
\usepackage{amssymb}
\usepackage{mathrsfs}

\def\thebibliograph#1#2{\section*{{\normalsize \bf #2}}\list
   {[\arabic{enumi}]}{\settowidth\labelwidth{[#1]}\leftmargin\labelwidth
     \advance\leftmargin\labelsep
     \usecounter{enumi}}
     \def\newblock{\hskip .11em plus .33em minus -.07em}
     \sloppy
     \sfcode`\.=1000\relax}

\newtheorem{theorem}{Theorem}

\newtheorem{definition}{Definition}

\newtheorem{lemma}{Lemma}

\input{tcilatex}
\begin{document}

\title{Real interpolation with variable exponent}
\author{Douadi Drihem}
\date{\today }
\maketitle

\begin{abstract}
We present the real interpolation with variable exponent and we prove the
basic properties in analogy to the classical real interpolation. More
precisely, we prove that under some additional conditions, this method\ can
be reduced to the case of fixed exponent. An application, we give the real
interpolation of variable Besov and Lorentz spaces as introduced recently in
Almeida and H\"{a}st\"{o} (J. Funct. Anal. 258 (5) 1628--2655, 2010) and L.
Ephremidze et al. (Fract. Calc. Appl. Anal. 11 (4) (2008), 407--420).\vskip%
5pt


\end{abstract}

\section{Introduction}

It is well known that real interpolation play an important role in several
different areas, especially for modern analysis and its theory started early
in 1960's by J-L. Lions and J. Peetre. There are two ways for introducing
the real interpolation method. The first is the $K$-method and the second is
the $J$-method. But the spaces generated by the $K$-and $J$-methods are the
same. For general literature on real interpolation we refer to \cite{BS88}, 
\cite{BL76}, \cite{T78} and references therein.

In recent years, there has been growing interest in generalizing classical
spaces such as Lebesgue spaces, Sobolev spaces, Besov spaces,
Triebel-Lizorkin spaces to the case with either variable integrability or
variable smoothness. The motivation for the increasing interest in such
spaces comes not only from theoretical purposes, but also from applications
to fluid dynamics, image restoration and PDE\ with non-standard growth
conditions.

From these in this paper we present a variable version of real
interpolation. First we study the variable version of $K$-method, where we
present some equivalent norms for the space generated by this method and we
prove their basic properties in analogy to the fixed exponent. Secondly, we
present the same analysis for the variable version of $J$-method and we
prove the first main statement of this paper. That is, under some additional
conditions the spaces generated by the $K$-and $J$-methods are the same.
Since the reiteration theorem is one of the most important general results
in interpolation theory, we will give its proof. Finally, we study the real
interpolation of variable exponent Besov and Lorentz spaces. Almost all of
the material we present is due to \cite{BS88}, and \cite{BL76}. Allowing the
exponent is vary from point to point will raise extra difficulties which, in
general, are overcome by imposing regularity assumptions on this exponent.

\section{Preliminaries}

As usual, we denote by $\mathbb{R}^{n}$ the $n$-dimensional real Euclidean
space, $\mathbb{N}$ the collection of all natural numbers and $\mathbb{N}%
_{0}=\mathbb{N}\cup \{0\}$. The letter $\mathbb{Z}$ stands for the set of
all integer numbers.\ The expression $f\lesssim g$ means that $f\leq c\,g$
for some independent constant $c$ (and non-negative functions $f$ and $g$),
and $f\approx g$ means $f\lesssim g\lesssim f$. \vskip5pt

By $c$ we denote generic positive constants, which may have different values
at different occurrences. Although the exact values of the constants are
usually irrelevant for our purposes, sometimes we emphasize their dependence
on certain parameters (e.g. $c(p)$ means that $c$ depends on $p$, etc.).
Further notation will be properly introduced whenever needed.

The variable exponents that we consider are always measurable functions $p$
on $\mathbb{R}$ with range in $[1,\infty \lbrack $. We denote the set of
such functions by $\mathcal{P}({\mathbb{R}})$. We use the standard notation $%
p^{-}:=\underset{x\in \mathbb{R}^{n}}{\text{ess-inf}}$ $p(x)$,$\quad p^{+}:=%
\underset{x\in \mathbb{R}^{n}}{\text{ess-sup }}p(x)$.

The variable exponent modular is defined by $\varrho _{p(\cdot )}(f):=\int_{%
\mathbb{R}^{n}}\varrho _{p(x)}(\left\vert f(x)\right\vert )dx$, where $%
\varrho _{p}(t)=t^{p}$. The variable exponent Lebesgue space $L^{p(\cdot )}$%
\ consists of measurable functions $f$ on $\mathbb{R}^{n}$ such that $%
\varrho _{p(\cdot )}(\lambda f)<\infty $ for some $\lambda >0$. We define
the Luxemburg (quasi)-norm on this space by the formula $\left\Vert
f\right\Vert _{p(\cdot )}:=\inf \Big\{\lambda >0:\varrho _{p(\cdot )}\Big(%
\frac{f}{\lambda }\Big)\leq 1\Big\}$. A useful property is that $\left\Vert
f\right\Vert _{p(\cdot )}\leq 1$ if and only if $\varrho _{p(\cdot )}(f)\leq
1$, see \cite{DHHR}, Lemma 3.2.4.

We say that $g:\mathbb{R}^{n}\rightarrow \mathbb{R}$ is \textit{locally }log%
\textit{-H\"{o}lder continuous}, abbreviated $g\in C_{\text{loc}}^{\log }(%
\mathbb{R}^{n})$, if there exists $c_{\log }(g)>0$ such that%
\begin{equation}
\left\vert g(x)-g(y)\right\vert \leq \frac{c_{\log }(g)}{\log
(e+1/\left\vert x-y\right\vert )}  \label{lo-log-Holder}
\end{equation}%
for all $x,y\in \mathbb{R}^{n}$. If 
\begin{equation*}
|g(x)-g(0)|\leq \frac{c_{\log }}{\ln (e+1/|x|)}
\end{equation*}%
for all $x\in \mathbb{R}^{n}$, then we say that $g$ is \emph{$\log $-H\"{o}%
lder continuous at the origin} (or has a \emph{$\log $ decay at the origin}%
). We say that $g$ satisfies the log\textit{-H\"{o}lder decay condition}, if
there exists $g_{\infty }\in \mathbb{R}$ and a constant $c_{\log }>0$ such
that%
\begin{equation*}
\left\vert g(x)-g_{\infty }\right\vert \leq \frac{c_{\log }}{\log
(e+\left\vert x\right\vert )}
\end{equation*}%
for all $x\in \mathbb{R}^{n}$. We say that $g$ is \textit{globally}-log%
\textit{-H\"{o}lder continuous}, abbreviated $g\in C^{\log }$, if it is%
\textit{\ }locally log-H\"{o}lder continuous and satisfies the log-H\"{o}%
lder decay\textit{\ }condition.\textit{\ }The constants $c_{\log }(g)$ and $%
c_{\log }$ are called the \textit{locally }log\textit{-H\"{o}lder constant }%
and the log\textit{-H\"{o}lder decay constant}, respectively\textit{.} We
note that all functions $g\in C_{\text{loc}}^{\log }(\mathbb{R}^{n})$ always
belong to $L^{\infty }$.\vskip5pt

We refer to the recent monograph $\mathrm{\cite{CF13}}$ for further
properties, historical remarks and references on variable exponent spaces.

\subsection{Technical lemmas}

In this subsection we present some results which are useful for us. The
following lemma is from \cite{DHHR}.

\begin{lemma}
\label{estimate -modular}Let$\ A\subset \mathbb{R}^{n}\ $and$\ p\in \mathcal{%
P}(\mathbb{R}^{n})$\ with$\ p^{-}<\infty $. If $\varrho _{p\left( \cdot
\right) }(f\chi _{A})>0$ or $p^{+}<\infty $, then%
\begin{equation*}
\min \left\{ \left( \varrho _{p\left( \cdot \right) }(f\chi _{A})\right) ^{%
\frac{1}{p^{-}}},\left( \varrho _{p\left( \cdot \right) }(f\chi _{A})\right)
^{\frac{1}{p^{+}}}\right\} \leq \left\Vert f\right\Vert _{p\left( \cdot
\right) }\leq \max \left\{ \left( \varrho _{p\left( \cdot \right) }(f\chi
_{A})\right) ^{\frac{1}{p^{-}}},\left( \varrho _{p\left( \cdot \right)
}(f\chi _{A})\right) ^{\frac{1}{p^{+}}}\right\} .
\end{equation*}
\end{lemma}

The next lemma is a Hardy-type inequality which is easy to prove.

\begin{lemma}
\label{Hardy-inequality}\textit{Let }$0<a<1$ \textit{and }$0<q\leq \infty $%
\textit{. Let }$\left\{ \varepsilon _{k}\right\} _{k}$\textit{\ be a
sequences of positive real numbers and denote} $\delta _{k}=\sum_{j=-\infty
}^{\infty }a^{\left\vert k-j\right\vert }\varepsilon _{j}$.\textit{\ }Then
there exists constant $c>0\ $\textit{depending only on }$a$\textit{\ and }$q$
such that%
\begin{equation*}
\Big(\sum\limits_{k=-\infty }^{\infty }\delta _{k}^{q}\Big)^{1/q}\leq c\text{
}\Big(\sum\limits_{k=-\infty }^{\infty }\varepsilon _{k}^{q}\Big)^{1/q}.
\end{equation*}
\end{lemma}

We will make use of the following statement, see \cite{DHHMS}, Lemma 3.3 for 
$w:=1$.

\begin{lemma}
\label{DHHR-estimate}Let $Q=(a,b)\subset \mathbb{R}$ with $0<a<b<\infty $.
Let $p\in \mathcal{P}(\mathbb{R})$ \textit{be log-H\"{o}lder continuous at
the origin }and $w:\mathbb{R}\rightarrow \mathbb{R}^{+}$ be a weight
function. Then for every $m>0$ there exists $\gamma =e^{-4mc_{\log
}(1/p)}\in \left( 0,1\right) $ such that%
\begin{eqnarray*}
&&\Big(\frac{\gamma }{w(Q)}\int_{Q}\left\vert f(y)\right\vert w(y)dy\Big)%
^{p\left( x\right) } \\
&\leq &\max \big(1,\left( w(Q)\right) ^{1-\frac{p\left( x\right) }{p^{-}}}%
\big)\frac{1}{w(Q)}\int_{Q}\left\vert f(y)\right\vert ^{p\left( y,0\right)
}w(y)dy \\
&&+\omega (m,b)\Big(\frac{1}{w(Q)}\int_{Q}g(x,y)w(y)dy\Big)
\end{eqnarray*}%
hold if $0<w(Q)<\infty $, all $x\in Q\subset \mathbb{R}$ and all $f\in
L^{p\left( \cdot \right) }(w)+L^{\infty }$\ with $\left\Vert fw^{1/p\left(
\cdot \right) }\right\Vert _{p\left( \cdot \right) }+\left\Vert f\right\Vert
_{\infty }\leq 1$, where%
\begin{equation*}
\omega (m,b)=\min \left( b^{m},1\right) \text{, }p\left( y,0\right) =p\left(
y\right) \text{ \ \ and \ \ }g(x,y)=(e+\frac{1}{x})^{-m}+(e+\frac{1}{y})^{-m}
\end{equation*}%
or%
\begin{equation*}
\omega (m,b)=\min \left( b^{m},1\right) \text{, }p\left( y,0\right) =p\left(
0\right) \text{ \ \ and \ \ }g(x,y)=(e+\frac{1}{x})^{-m}\chi _{\{x:p(x){<}%
p(0)\}}(x).
\end{equation*}%
\textit{In addition we have the same estimate}, where%
\begin{equation*}
\omega (m,b)=1\text{, }p\left( y,0\right) =p_{\infty }\text{ \ \ and \ \ }%
g(x,y)=(e+x)^{-m}\chi _{\{x:p(x){<}p_{\infty }\}}(x),
\end{equation*}%
if $p\in \mathcal{P}(\mathbb{R})$ satisfies the log\textit{-H\"{o}lder decay
condition, where }$\gamma =e^{-4mc_{\log }}$\textit{.}
\end{lemma}

The proof of this lemma\ is given in \cite{D4}. Notice that in the proof of
this theorem we need only that 
\begin{equation*}
\int_{Q}\left\vert f(y)\right\vert ^{p\left( y\right) }w(y)dy\leq 1
\end{equation*}%
and/or $\left\Vert f\right\Vert _{\infty }\leq 1$.

The next lemma is the continuous version of Hardy-type inequality, see \cite%
{D.S2007}.

\begin{lemma}
\label{lq-inequality}\textit{Let }$s>0$\textit{. Let }$q\in \mathcal{P}(%
\mathbb{R})$\textit{\ }be \emph{$\log $}-H\"{o}lder continuous both at the
origin\ and at infinity with $1\leq q^{-}\leq q^{+}<\infty $\textit{. Let }$%
\left\{ \varepsilon _{t}\right\} _{t}$\textit{\ be a sequence of positive
measurable functions.} Let%
\begin{equation*}
\eta _{t}=t^{s}\int_{t}^{\infty }\tau ^{-s}\varepsilon _{\tau }\frac{d\tau }{%
\tau }\quad \text{and\quad }\delta _{t}=t^{-s}\int_{0}^{t}\tau
^{s}\varepsilon _{\tau }\frac{d\tau }{\tau }.
\end{equation*}%
Then there exists constant $c>0\ $\textit{depending only on }$s$, $q^{-}$%
\textit{, c}$_{\log }(q)$ \textit{and }$q^{+}$ such that%
\begin{equation*}
\left\Vert \eta _{t}\right\Vert _{L^{q(\cdot )}((0,\infty ),\frac{dt}{t}%
)}+\left\Vert \delta _{t}\right\Vert _{L^{q(\cdot )}((0,\infty ),\frac{dt}{t}%
)}\lesssim \left\Vert \varepsilon _{t}\right\Vert _{L^{q(\cdot )}((0,\infty
),\frac{dt}{t})}.
\end{equation*}
\end{lemma}

\section{The K-Method}

The fundamental notion of real interpolation is the $K$-functional, where it
is due to J. Peetre.

\begin{definition}
\label{Compatible}Let $A_{0}$ and $A_{1}$\ be Banach spaces over $\mathbb{K=R%
}$ or $\mathbb{C}$. We shall say that $A_{0}$ and $A_{1}$ are compatible if
there is a Hausdorff topological vector space $Z$ such that 
\begin{equation*}
A_{0},A_{1}\hookrightarrow Z,
\end{equation*}%
with continuous embeddings.
\end{definition}

Let $A_{0}$ and $A_{1}$ be compatible. We will say that $(A_{0},A_{1})$ is a
compatible couple. Then we can form their sum $A_{0}+A_{1}$ and their
intersection $A_{0}\cap A_{1}$. The sum consists of all $f\in Z$ such that
we can write

\begin{equation*}
f=f_{0}+f_{1}
\end{equation*}%
\ for some $f_{0}\in A_{0}$\ and $f_{1}\in A_{1}$. Then $A_{0}+A_{1}$ is a
Banach space with norm defined by\ 
\begin{equation*}
\left\Vert f\right\Vert _{A_{0}+A_{1}}=\inf_{f=f_{0}+f_{1}}\left( \left\Vert
f_{0}\right\Vert _{A_{0}}+\left\Vert f_{1}\right\Vert _{A_{1}}\right) .
\end{equation*}%
$A_{0}\cap A_{1}$ is a Banach space with norm defined by\ 
\begin{equation*}
\left\Vert f\right\Vert _{A_{0}\cap A_{1}}=\max_{f=f_{0}+f_{1}}\left(
\left\Vert f_{0}\right\Vert _{A_{0}},\left\Vert f_{1}\right\Vert
_{A_{1}}\right) .
\end{equation*}

Let $(A_{0},A_{1})$ be a compatible couple. With $t>0$\ fixed, put%
\begin{equation*}
K(t,f,A_{0},A_{1})=\inf_{f=f_{0}+f_{1}}\left( \left\Vert f_{0}\right\Vert
_{A_{0}}+t\left\Vert f_{1}\right\Vert _{A_{1}}\right) ,\quad f\in
A_{0}+A_{1},
\end{equation*}

is the $K$-functional. For any $f\in A_{0}+A_{1}$, $K(t,f,A_{0},A_{1})$ is a
positive, increasing and concave function of $t$. In particular%
\begin{equation}
K(t,f,A_{0},A_{1})\leq \max (1,\frac{t}{s})K(s,f,A_{0},A_{1}).
\label{K-property}
\end{equation}%
If there is no danger of confusion, we shall write $%
K(t,f)=K(t,f,A_{0},A_{1}) $.

\begin{definition}
Let $\theta \in (0,1)$ and $q\in \mathcal{P}(\mathbb{R})$. Let $%
(A_{0},A_{1}) $\ be a compatible couple. The space $(A_{0},A_{1})_{\theta
,q(\cdot )}$ consists of all $f$ in $A_{0}+A_{1}$ for which the functional%
\begin{equation*}
\left\Vert f\right\Vert _{(A_{0},A_{1})_{\theta ,q(\cdot )}}=\left\Vert
t^{-\theta }K(t,f)\right\Vert _{L^{q(\cdot )}((0,\infty ),\frac{dt}{t})}
\end{equation*}%
is finite.
\end{definition}

\begin{definition}
Let $\theta \in \lbrack 0,1]$. Let $(A_{0},A_{1})$\ be a compatible couple.
The space $(A_{0},A_{1})_{\theta ,\infty }$ consists\ of\ all\ $f$\ in $%
A_{0}+A_{1}$ for which the functional%
\begin{equation*}
\left\Vert f\right\Vert _{(A_{0},A_{1})_{\theta ,\infty
}}=\sup_{t>0}t^{-\theta }K(t,f)
\end{equation*}%
is finite.
\end{definition}

In the next lemma we prove that the first definition can be given in
discrete\ version, where we need additional assumptions on $q$.

\begin{lemma}
\label{K-property2}Let $(A_{0},A_{1})$\ be a compatible couple and $f\in
A_{0}+A_{1}$.\ Let $\theta \in (0,1)$, $f\in (A_{0},A_{1})_{\theta ,q(\cdot
)}$\ and we put $\alpha _{v}=K(2^{v},f)$, $v\in \mathbb{Z}$\textit{. Let }$%
q\in \mathcal{P}(\mathbb{R})$ \textit{be log-H\"{o}lder continuous both at
the origin and at infinity. Then}%
\begin{equation*}
\left\Vert f\right\Vert _{(A_{0},A_{1})_{\theta ,q(\cdot )}}\approx \Big(%
\sum_{v=-\infty }^{0}2^{-v\theta q(0)}\alpha _{v}^{q(0)}\Big)^{\frac{1}{q(0)}%
}+\Big(\sum_{v=1}^{\infty }2^{-v\theta q_{\infty }}\alpha _{v}^{q_{\infty }}%
\Big)^{\frac{1}{q_{\infty }}}.
\end{equation*}%
Moreover,%
\begin{equation*}
\left\Vert f\right\Vert _{(A_{0},A_{1})_{\theta ,q(\cdot )}}\approx \Big(%
\int_{0}^{1}t^{-\theta q(0)}K(t,f)^{q(0)}\frac{dt}{t}\Big)^{\frac{1}{q(0)}}+%
\Big(\int_{1}^{\infty }t^{-\theta q_{\infty }}K(t,f)^{q_{\infty }}\frac{dt}{t%
}\Big)^{\frac{1}{q_{\infty }}}.
\end{equation*}
\end{lemma}

\textbf{Proof. }We will do the proof in two steps and we need only to prove
the first statement.

\textit{Step 1.} Let us prove that%
\begin{equation}
S=\Big(\sum_{v=-\infty }^{0}2^{-v\theta q(0)}\alpha _{v}^{q(0)}\Big)^{\frac{1%
}{q(0)}}+\Big(\sum_{v=1}^{\infty }2^{-v\theta q_{\infty }}\alpha
_{v}^{q_{\infty }}\Big)^{\frac{1}{q_{\infty }}}\lesssim \left\Vert
f\right\Vert _{(A_{0},A_{1})_{\theta ,q(\cdot )}}.  \label{est-step1}
\end{equation}%
By scaling argument, we need only to prove that%
\begin{equation*}
\sum_{v=1}^{\infty }2^{-v\theta q_{\infty }}\big(\alpha _{v}\big)^{q_{\infty
}}\lesssim 1\quad \text{and}\quad \sum_{v=-\infty }^{0}2^{-v\theta q(0)}\big(%
\alpha _{v}\big)^{q(0)}\lesssim 1
\end{equation*}%
for any $f\in (A_{0},A_{1})_{\theta ,q(\cdot )}$ with $\left\Vert
f\right\Vert _{(A_{0},A_{1})_{\theta ,q(\cdot )}}\leq 1$. To prove the first
estimate we need to prove that 
\begin{equation*}
2^{-v\theta q_{\infty }}\big(\alpha _{v}\big)^{q_{\infty }}\leq
\int_{2^{v-1}}^{2^{v}}\big(t^{-\theta }K(t,f)\big)^{q\left( t\right) }\frac{%
dt}{t}+2^{-v}=\delta
\end{equation*}%
for any $v\in \mathbb{N}$. This claim can be reformulated as showing that%
\begin{equation*}
\big(\delta ^{-\frac{1}{q_{\infty }}}2^{-v\theta }\alpha _{v}\big)%
^{q_{\infty }}=\Big(\frac{1}{\log 2}\int_{2^{v-1}}^{2^{v}}\delta ^{-\frac{1}{%
q_{\infty }}}2^{-v\theta }\alpha _{v}\frac{d\tau }{\tau }\Big)^{q_{\infty
}}\lesssim 1.
\end{equation*}%
Using the property $\mathrm{\eqref{K-property}}$, we find that%
\begin{equation*}
\int_{2^{v-1}}^{2^{v}}\delta ^{-\frac{1}{q_{\infty }}}2^{-v\theta }\alpha
_{v}\frac{d\tau }{\tau }\lesssim \int_{2^{v-1}}^{2^{v}}\delta ^{-\frac{1}{%
q_{\infty }}}\tau ^{-\theta }K(\tau ,f)\frac{d\tau }{\tau }.
\end{equation*}%
By Lemma \ref{DHHR-estimate} the last expression with power $q\left(
t\right) $ is bounded by%
\begin{equation*}
c\int_{2^{v-1}}^{2^{v}}\delta ^{-\frac{q\left( \tau \right) }{q_{\infty }}}%
\big(\tau ^{-\theta }K(\tau ,f)\big)^{q\left( \tau \right) }\frac{d\tau }{%
\tau }+1
\end{equation*}%
for any $t\in \lbrack 2^{v-1},2^{v}]$. Since $q$ is \emph{$\log $-}H\"{o}%
lder continuous at infinity, we find that%
\begin{equation}
\delta ^{-\frac{q\left( \tau \right) }{q_{\infty }}}\approx \delta ^{-1},%
\text{ \ \ }\tau \in \lbrack 2^{v-1},2^{v}],\quad v\in \mathbb{N}.
\label{log-at- zero}
\end{equation}%
Therefore, from the definition of $\delta $, we find that%
\begin{equation*}
\int_{2^{v-1}}^{2^{v}}\delta ^{-\frac{q\left( \tau \right) }{q_{\infty }}}%
\big(\tau ^{-\theta }K(\tau ,f)\big)^{q\left( \tau \right) }\frac{d\tau }{%
\tau }\lesssim 1.
\end{equation*}%
Now, let us prove the second estimate. We need to show that 
\begin{equation*}
2^{-v\theta q(0)}\big(\alpha _{v}\big)^{q(0)}\lesssim \int_{2^{v}}^{2^{v+1}}%
\big(t^{-\theta }K(t,f)\big)^{q\left( t\right) }\frac{dt}{t}+2^{v}=\delta
\end{equation*}%
for any $v\leq 0$. This claim can be reformulated as showing that%
\begin{equation*}
\big(\delta ^{-\frac{1}{q(0)}}2^{-v\theta }\alpha _{v}\big)^{q(0)}=\Big(%
\frac{1}{\log 2}\int_{2^{v}}^{2^{v+1}}\delta ^{-\frac{1}{q(0)}}2^{-v\theta
}\alpha _{v}\frac{d\tau }{\tau }\Big)^{q(0)}\lesssim 1.
\end{equation*}%
The property $\mathrm{\eqref{K-property}}$, gives\ that%
\begin{equation*}
\int_{2^{v}}^{2^{v+1}}\delta ^{-\frac{1}{q(0)}}2^{-v\theta }\alpha _{v}\frac{%
d\tau }{\tau }\lesssim \int_{2^{v}}^{2^{v+1}}\delta ^{-\frac{1}{q(0)}}\tau
^{-\theta }K(\tau ,f)\frac{d\tau }{\tau }.
\end{equation*}%
Again by Lemma \ref{DHHR-estimate},%
\begin{equation*}
\Big(\int_{2^{v}}^{2^{v+1}}\delta ^{-\frac{1}{q(0)}}\tau ^{-\theta }K(\tau
,f)\frac{d\tau }{\tau }\Big)^{q\left( t\right) }\lesssim
\int_{2^{v}}^{2^{v+1}}\delta ^{-\frac{q\left( \tau \right) }{q(0)}}\big(\tau
^{-\theta }K(\tau ,f)\big)^{q\left( \tau \right) }\frac{d\tau }{\tau }+1.
\end{equation*}%
for any $t\in \lbrack 2^{v},2^{v+1}]$\ and any $v\leq 0$. We use the
logarithmic decay\ condition at origin\ of $q$ to show that%
\begin{equation*}
\delta ^{-\frac{q\left( \tau \right) }{q(0)}}\approx \delta ^{-1},\text{ \ \ 
}\tau \in \lbrack 2^{v},2^{v+1}],\quad v\leq 0.
\end{equation*}%
Therefore, from the definition of $\delta $, we find that%
\begin{equation*}
\int_{2^{v}}^{2^{v+1}}\delta ^{-\frac{q\left( \tau \right) }{q(0)}}\big(\tau
^{-\theta }K(\tau ,f)\big)^{q\left( \tau \right) }\frac{d\tau }{\tau }%
\lesssim 1
\end{equation*}%
for any $v\leq 0$. Hence, we proved $\mathrm{\eqref{est-step1}}$.

\textit{Step 2.} Let us prove that%
\begin{equation*}
\left\Vert f\right\Vert _{(A_{0},A_{1})_{\theta ,q(\cdot )}}\lesssim S.
\end{equation*}%
This claim can be reformulated as showing that%
\begin{equation*}
\dint_{0}^{\infty }\big(t^{-\theta }K(t,\frac{f}{S})\big)^{q\left( t\right) }%
\frac{dt}{t}\lesssim 1.
\end{equation*}%
Now our estimate clearly follows from the inequalities%
\begin{equation*}
\int_{2^{v-1}}^{2^{v}}\big(t^{-\theta }K(t,\frac{f}{S})\big)^{q\left(
t\right) }\frac{dt}{t}\lesssim 2^{-v\theta q_{\infty }}\big(\frac{\alpha _{v}%
}{S}\big)^{q_{\infty }}+2^{-v}=\delta
\end{equation*}%
for any $v\in \mathbb{N}$ and 
\begin{equation*}
\int_{2^{v}}^{2^{v+1}}\big(t^{-\theta }K(t,\frac{f}{S})\big)^{q\left(
t\right) }\frac{dt}{t}\lesssim 2^{-v\theta q(0)}\big(\frac{\alpha _{v}}{S}%
\big)^{q(0)}+2^{v}
\end{equation*}%
for any $v\leq 0$. The first claim can be reformulated as showing\ that%
\begin{equation*}
\int_{2^{v-1}}^{2^{v}}\big(\delta ^{-\frac{1}{q(t)}}t^{-\theta }K(t,\frac{f}{%
S})\big)^{q\left( t\right) }\frac{dt}{t}\lesssim 1.
\end{equation*}%
We need only to show that%
\begin{equation*}
\delta ^{-\frac{1}{q(t)}}t^{-\theta }K(t,\frac{f}{S})\lesssim 1
\end{equation*}%
for any $v\in \mathbb{N}$ and any $t\in \lbrack 2^{v-1},2^{v}]$. From $%
\mathrm{\eqref{K-property}}$, the left-hand side is bounded by%
\begin{equation*}
\delta ^{-\frac{1}{q(t)}}2^{-\theta v}K(2^{v},\frac{f}{S}),
\end{equation*}%
and from $\mathrm{\eqref{log-at- zero}}$ we find that 
\begin{equation*}
\delta ^{-\frac{1}{q(t)}}2^{-\theta v}K(2^{v},\frac{f}{S})\lesssim \delta ^{-%
\frac{1}{q_{\infty }}}2^{-\theta v}K(2^{v},\frac{f}{S})\leq 1.
\end{equation*}%
for any $v\in \mathbb{N}$. Similarly we estimate the second claim. Hence the
lemma is proved.\quad $\square $

Let $(A_{0},A_{1})$\ be a compatible couple and\ $f\in A_{0}+A_{1}$.\ Let $%
\theta \in \lbrack 0,1]$, $f\in (A_{0},A_{1})_{\theta ,q(\cdot )}$\ and we
put $\alpha _{v}=K(2^{v},f)$, $v\in \mathbb{Z}$\textit{. }Then we have%
\begin{equation*}
\left\Vert f\right\Vert _{(A_{0},A_{1})_{\theta ,\infty }}\approx \sup_{v\in 
\mathbb{Z}}2^{v\theta }\alpha _{v}.
\end{equation*}

We present\ some\ important properties of the spaces $(A_{0},A_{1})_{\theta
,q(\cdot )}$.

\begin{theorem}
\label{embedding}Let\ $\theta \in (0,1)$\ and\ $q\in \mathcal{P}(\mathbb{R})$%
. Let $(A_{0},A_{1})$\ be a compatible couple of Banach spaces. Then $%
(A_{0},A_{1})_{\theta ,q(\cdot )}$ is Banach space and 
\begin{equation*}
K(s,f,A_{0},A_{1})\leq \gamma _{\theta ,q^{+}}s^{\theta }\left\Vert
f\right\Vert _{(A_{0},A_{1})_{\theta ,q(\cdot )}}
\end{equation*}%
for any $s>0$. Moreover we have%
\begin{equation*}
A_{0}\cap A_{1}\hookrightarrow \left\Vert f\right\Vert
_{(A_{0},A_{1})_{\theta ,q(\cdot )}}\hookrightarrow A_{0}+A_{1}.
\end{equation*}
\end{theorem}

\textbf{Proof.}\ Let $\{f_{n}\}_{n}$ be a sequence in $A_{0}+A_{1}$ such that%
\begin{equation*}
\sum_{n=1}^{\infty }\left\Vert f_{n}\right\Vert _{(A_{0},A_{1})_{\theta
,q(\cdot )}}<\infty .
\end{equation*}%
Since $L^{q(\cdot )}((0,\infty ),\frac{dt}{t})$ is a Banach space, the
series $\sum_{n=1}^{\infty }t^{-\theta }K(t,f_{n})$\ converges in $%
L^{q(\cdot )}((0,\infty ),\frac{dt}{t})$ then we get%
\begin{equation*}
\sum_{n=1}^{\infty }t^{-\theta }K(t,f_{n})<\infty
\end{equation*}%
for all $t>0$. Since $A_{0}+A_{1}$\ is a Banach space, then%
\begin{equation*}
t^{-\theta }K(t,\sum_{n=1}^{\infty }f_{n})\leq t^{-\theta
}\sum_{n=1}^{\infty }K(t,f_{n})
\end{equation*}%
for all $t>0$. Applying the $L^{q(\cdot )}((0,\infty ),\frac{dt}{t})$-norm
to each side, we obtain 
\begin{equation*}
\big\|\sum_{n=1}^{\infty }f_{n}\big\|_{(A_{0},A_{1})_{\theta ,q(\cdot
)}}\leq \sum_{n=1}^{\infty }\left\Vert f_{n}\right\Vert
_{(A_{0},A_{1})_{\theta ,q(\cdot )}}<\infty ,
\end{equation*}%
which ensure that $(A_{0},A_{1})_{\theta ,q(\cdot )}$ is Banach space. By
the property $\mathrm{\eqref{K-property}}$ we find that%
\begin{equation*}
\min (1,\frac{t}{s})K(s,f)\leq K(t,f),\quad s,t>0.
\end{equation*}%
Therefore,%
\begin{equation*}
\big\|t^{-\theta }\min (1,\frac{t}{s})\big\|_{L^{q(\cdot )}((0,\infty ),%
\frac{dt}{t})}K(s,f)\leq \left\Vert f\right\Vert _{(A_{0},A_{1})_{\theta
,q(\cdot )}}.
\end{equation*}%
Let us prove that%
\begin{equation}
\big\|t^{-\theta }\min (1,\frac{t}{s})\big\|_{L^{q(\cdot )}((0,\infty ),%
\frac{dt}{t})}\gtrsim s^{-\theta }.  \label{est-t-theta}
\end{equation}%
We have%
\begin{equation*}
\big\|t^{-\theta }\min (1,\frac{t}{s})\big\|_{L^{q(\cdot )}((0,\infty ),%
\frac{dt}{t})}\geq s^{-\theta }\big\|\big(\frac{t}{s}\big)^{1-\theta }\big\|%
_{L^{q(\cdot )}((0,s),\frac{dt}{t})},
\end{equation*}%
and%
\begin{equation*}
\int_{0}^{s}\big(\frac{t}{s}\big)^{(1-\theta )q(t)}\frac{dt}{t}\geq
\int_{0}^{s}\big(\frac{t}{s}\big)^{(1-\theta )q^{+}}\frac{dt}{t}=\frac{1}{%
(1-\theta )q^{+}}.
\end{equation*}%
From Lemma \ref{estimate -modular}, we find our claim $\mathrm{%
\eqref{est-t-theta}}$. Therefore,%
\begin{equation*}
K(s,f)\leq \gamma _{\theta ,q^{+}}s^{\theta }\left\Vert f\right\Vert
_{(A_{0},A_{1})_{\theta ,q(\cdot )}}
\end{equation*}%
for any $s>0$. Taking $s=1$, we obtain 
\begin{equation*}
\left\Vert f\right\Vert _{A_{0}+A_{1}}\lesssim \left\Vert f\right\Vert
_{(A_{0},A_{1})_{\theta ,q(\cdot )}}.
\end{equation*}%
Now since%
\begin{equation*}
K(t,f)\leq \min (1,t)\left\Vert f\right\Vert _{A_{0}\cap A_{1}},
\end{equation*}%
we find that%
\begin{equation*}
\left\Vert f\right\Vert _{(A_{0},A_{1})_{\theta ,q(\cdot )}}\lesssim
\left\Vert f\right\Vert _{A_{0}\cap A_{1}}.
\end{equation*}%
\quad $\square $

\begin{definition}
Let $(A_{0},A_{1})$\ and $(B_{0},B_{1})$ be two compatible couples of Banach
spaces and let $T$ be a linear operator defined on $A_{0}+A_{1}$ and taking
values in $B_{0}+B_{1}$. $T$ is said be admissible\ with respect to the
couples\ $(A_{0},A_{1})$\ and $(B_{0},B_{1})$ if, for each $i=1,0$ the
restriction of $T$ to $A_{i}$ maps $A_{i}$ into $B_{i}$ and furthermore is a
bounded operator from $A_{i}$ into $B_{i}:$%
\begin{equation*}
\left\Vert Tf\right\Vert _{B_{i}}\leq \left\Vert T\right\Vert
_{L(A_{i},B_{i})}\left\Vert f\right\Vert _{A_{i}},\quad f\in A_{i}.
\end{equation*}
\end{definition}

Notice that every admissible operator $T$\ with respect to the couples\ $%
(A_{0},A_{1})$\ and $(B_{0},B_{1})$ is bounded from $A_{0}+A_{1}$ into $%
B_{0}+B_{1}$.

\begin{theorem}
\label{b-operators}Let\ $\theta \in (0,1)$\ and\ $q\in \mathcal{P}(\mathbb{R}%
)$. Let $(A_{0},A_{1})$\ and $(B_{0},B_{1})$ be two compatible couples of
Banach spaces and let $T$ be admissible\ with respect to the couples\ $%
(A_{0},A_{1})$\ and $(B_{0},B_{1})$. Then%
\begin{equation*}
T:(A_{0},A_{1})_{\theta ,q(\cdot )}\longrightarrow (B_{0},B_{1})_{\theta
,q(\cdot )}
\end{equation*}%
and%
\begin{equation*}
\left\Vert Tf\right\Vert _{(B_{0},B_{1})_{\theta ,q(\cdot )}}\leq \max
\left( \left\Vert T\right\Vert _{L(A_{0},B_{0})},\left\Vert T\right\Vert
_{L(A_{1},B_{1})}\right) \left\Vert f\right\Vert _{(A_{0},A_{1})_{\theta
,q(\cdot )}}
\end{equation*}%
for all $f\in (A_{0},A_{1})_{\theta ,q(\cdot )}$.
\end{theorem}

\textbf{Proof.} Suppose that $T:(A_{0},A_{1})\longrightarrow (B_{0},B_{1})$.
Then%
\begin{eqnarray*}
K(t,Tf,B_{0},B_{1}) &\leq &\left\Vert T\right\Vert _{L(A_{0},B_{0})}K\Big(%
\frac{\left\Vert T\right\Vert _{L(A_{1},B_{1})}t}{\left\Vert T\right\Vert
_{L(A_{0},B_{0})}},f,A_{0},A_{1}\Big) \\
&\leq &\max \big(\left\Vert T\right\Vert _{L(A_{0},B_{0})},\left\Vert
T\right\Vert _{L(A_{1},B_{1})}\big)K(t,f,A_{0},A_{1}),
\end{eqnarray*}%
by the property $\mathrm{\eqref{K-property}}$. Multiplying by $t^{-\theta }$
and then applying the $L^{q(\cdot )}((0,\infty ),\frac{dt}{t})$-norm to each
side we obtain the desired estimate.\quad $\square $

\begin{proposition}
Let $\theta \in (0,1)$. Let $(A_{0},A_{1})$\ be a compatible couples of
Banach spaces. \newline
(i) Let $q,r\in \mathcal{P}(\mathbb{R})$ with $1\leq q(\cdot )\leq r(\cdot
)<\infty $. Then%
\begin{equation*}
(A_{0},A_{1})_{\theta ,q(\cdot )}\hookrightarrow (A_{0},A_{1})_{\theta
,r(\cdot )}.
\end{equation*}%
and%
\begin{equation*}
(A_{0},A_{1})_{\theta ,q(\cdot )}\hookrightarrow (A_{0},A_{1})_{\theta
,\infty }
\end{equation*}%
(ii) \textit{Let }$q\in \mathcal{P}(\mathbb{R})$ \textit{be log-H\"{o}lder
continuous both at the origin and at infinity with $q(0)=q_{\infty }$. Then}%
\begin{equation*}
(A_{0},A_{1})_{\theta ,q(\cdot )}=(A_{1},A_{0})_{1-\theta ,q(\cdot )}.
\end{equation*}%
(iii) \textit{Let }$q,r\in \mathcal{P}(\mathbb{R})$ \textit{be log-H\"{o}%
lder continuous both at the origin and at infinity with }$q(0)=r(0)$ and $%
q_{\infty }=r_{\infty }$\textit{. Then}%
\begin{equation*}
(A_{0},A_{1})_{\theta ,q(\cdot )}=(A_{0},A_{1})_{\theta ,r(\cdot )}.
\end{equation*}%
(iv) If $A_{1}\hookrightarrow A_{0}$, then%
\begin{equation*}
(A_{0},A_{1})_{\theta _{1},q(\cdot )}\hookrightarrow (A_{0},A_{1})_{\theta
,q(\cdot )}\text{\quad if }0<\theta \leq \theta _{1}<1.
\end{equation*}%
(v) If $A_{0}=A_{1}$, with equal norm, then%
\begin{equation*}
(A_{0},A_{1})_{\theta ,q(\cdot )}=A_{0}.
\end{equation*}
\end{proposition}

\textbf{Proof.} We prove (i). From Theorem \ref{embedding}, we obtain%
\begin{equation*}
(A_{0},A_{1})_{\theta ,q(\cdot )}\hookrightarrow (A_{0},A_{1})_{\theta
,\infty }\quad \text{and}\quad K(s,f)\lesssim s^{\theta }\left\Vert
f\right\Vert _{(A_{0},A_{1})_{\theta ,q(\cdot )}}
\end{equation*}%
for any $f\in (A_{0},A_{1})_{\theta ,q(\cdot )}$, any $s>0$\ with $%
\left\Vert f\right\Vert _{(A_{0},A_{1})_{\theta ,q(\cdot )}}\neq 0$\ and
this implies that%
\begin{eqnarray*}
&&\int_{0}^{\infty }\big(t^{-\theta }K(t,\frac{f}{\left\Vert f\right\Vert
_{(A_{0},A_{1})_{\theta ,q(\cdot )}}})\big)^{r\left( t\right) }\frac{dt}{t}
\\
&\leq &\int_{0}^{\infty }\big(t^{-\theta }K(t,\frac{f}{\left\Vert
f\right\Vert _{(A_{0},A_{1})_{\theta ,q(\cdot )}}})\big)^{q\left( t\right) }%
\big(\sup_{t>0}t^{-\theta }K(t,\frac{f}{\left\Vert f\right\Vert
_{(A_{0},A_{1})_{\theta ,q(\cdot )}}})\big)^{r\left( t\right) -q\left(
t\right) }\frac{dt}{t} \\
&\lesssim &\int_{0}^{\infty }\big(t^{-\theta }K(t,\frac{f}{\left\Vert
f\right\Vert _{(A_{0},A_{1})_{\theta ,q(\cdot )}}})\big)^{q\left( t\right) }%
\frac{dt}{t}.
\end{eqnarray*}%
The last term is bounded since 
\begin{equation*}
\Big\|t^{-\theta }K(t,\frac{f}{\left\Vert f\right\Vert
_{(A_{0},A_{1})_{\theta ,q(\cdot )}}})\Big\|_{L^{q(\cdot )}((0,\infty ),%
\frac{dt}{t})}=1,
\end{equation*}%
and hence%
\begin{equation*}
\left\Vert f\right\Vert _{(A_{0},A_{1})_{\theta ,r(\cdot )}}\lesssim
\left\Vert f\right\Vert _{(A_{0},A_{1})_{\theta ,q(\cdot )}}.
\end{equation*}%
Hence the property (i) is proved. To prove (ii) we use Lemma \ref%
{K-property2} and the fact that 
\begin{equation*}
K(t,f,A_{0},A_{1})=tK(t^{-1},f,A_{1},A_{0}),\quad v>0,
\end{equation*}%
and $q(0)=q_{\infty }$. The property (iii) follows by Lemma \ref{K-property2}%
. Now if $A_{1}\hookrightarrow A_{0}$ we have $\left\Vert f\right\Vert
_{A_{0}}\leq c\left\Vert f\right\Vert _{A_{1}}$ for any $f\in A_{0}$ and 
\begin{equation*}
K(t,f)=\left\Vert f\right\Vert _{A_{0}},
\end{equation*}%
if $t>c$. Then%
\begin{equation*}
\left\Vert t^{-\theta }K(t,f)\right\Vert _{L^{q(\cdot )}((c,\infty ),\frac{dt%
}{t})}\lesssim \left\Vert f\right\Vert _{A_{0}},
\end{equation*}%
and 
\begin{equation*}
\left\Vert t^{-\theta }K(t,f)\right\Vert _{L^{q(\cdot )}((0,\infty ),\frac{dt%
}{t})}\lesssim \left\Vert t^{-\theta }K(t,f)\right\Vert _{L^{q(\cdot
)}((0,c),\frac{dt}{t})}+\left\Vert f\right\Vert _{A_{0}}.
\end{equation*}%
Using the fact that%
\begin{equation*}
\left\Vert f\right\Vert _{A_{0}}\lesssim \left\Vert t^{-\theta
_{1}}K(t,f)\right\Vert _{L^{q(\cdot )}((c,\infty ),\frac{dt}{t})},
\end{equation*}%
and $0<\theta \leq \theta _{1}<1$, we obtain%
\begin{equation*}
\left\Vert f\right\Vert _{(A_{0},A_{1})_{\theta ,q(\cdot )}}\lesssim
\left\Vert f\right\Vert _{(A_{0},A_{1})_{\theta _{1},q(\cdot )}}.
\end{equation*}%
So, the property (iv) is proved. Now the property (v) is immediate. The
proof is complete.\quad $\square $

\section{The J-Method}

Let $(A_{0},A_{1})$ be a compatible couple. With $t>0$\ fixed, put%
\begin{equation*}
J(t,f,A_{0},A_{1})=\inf_{f=f_{0}+f_{1}}\left( \left\Vert f_{0}\right\Vert
_{A_{0}},t\left\Vert f_{1}\right\Vert _{A_{1}}\right) ,\quad f\in A_{0}\cap
A_{1}.
\end{equation*}%
Notice that $J(t,f,A_{0},A_{1})$\ is an equivalent norm on $A_{0}\cap A_{1}$
for a given $t>0$. If there is no danger of confusion, we shall write $%
J(t,f)=J(t,f,A_{0},A_{1})$. For any $f\in A_{0}\cap A_{1}$, $J(t,f)$ is a
positive, increasing and convex \ function of $t$, such that%
\begin{equation}
J(t,f)\leq \max (1,\frac{t}{s})J(s,f),  \label{Est-J}
\end{equation}%
and%
\begin{equation}
K(t,f)\leq \min (1,\frac{t}{s})J(s,f).  \label{Est-J1}
\end{equation}%
Now we define the interpolation space constructed by the $J$-method.

\begin{definition}
Let $\theta \in (0,1)$ and $q\in \mathcal{P}(\mathbb{R})$. Let $%
(A_{0},A_{1}) $ be a compatible couple. The space $(A_{0},A_{1})_{\theta
,q(\cdot ),J}$ consists of all $f$ in $A_{0}+A_{1}$ that are representable
in the form%
\begin{equation}
f=\int_{0}^{\infty }u(t)\frac{dt}{t}  \label{rep}
\end{equation}%
where $u(t)$ is measurable with values in $A_{0}\cap A_{1}$\ and%
\begin{equation*}
\left\Vert f\right\Vert _{(A_{0},A_{1})_{\theta ,q(\cdot ),J}}=\inf
\left\Vert t^{-\theta }J(t,u(t))\right\Vert _{L^{q(\cdot )}((0,\infty ),%
\frac{dt}{t})}<\infty ,
\end{equation*}%
where the infimum is taken over all $u$ such that $\mathrm{\eqref{rep}}$
holds.
\end{definition}

\begin{definition}
Let $\theta \in (0,1)$. The space $(A_{0},A_{1})_{\theta ,\infty ,J}$
consists of all $f$ in $A_{0}+A_{1}$ that are representable in the form $%
\mathrm{\eqref{rep}}$, where $u(t)$ is measurable with values in $A_{0}\cap
A_{1}$\ and%
\begin{equation*}
\left\Vert f\right\Vert _{(A_{0},A_{1})_{\theta ,\infty ,J}}=\inf
\sup_{t>0}t^{-\theta }J(t,u(t))<\infty ,
\end{equation*}%
where the infimum is taken over all $u$\ such that\ $\mathrm{\eqref{rep}}$
holds.
\end{definition}

\begin{lemma}
\label{J-property2}Let $(A_{0},A_{1})$\ be a compatible couple and $f\in
A_{0}+A_{1}$.\ Let $\theta \in (0,1)$\textit{\ and }$q\in \mathcal{P}(%
\mathbb{R})$ \textit{be log-H\"{o}lder continuous both at the origin and at
infinity. Then }$f\in (A_{0},A_{1})_{\theta ,q(\cdot ),J}$ if and only if
there exist $u_{v}\in A_{0}\cap A_{1}$, $v\in \mathbb{Z}$, with 
\begin{equation}
f=\sum_{v=-\infty }^{\infty }u_{v}\quad \text{convergence in }A_{0}\cap
A_{1},  \label{f-representation}
\end{equation}%
and such that 
\begin{equation*}
\left\Vert \left( J(2^{v},u_{v})\right) _{v}\right\Vert _{\lambda ^{\theta
,q(0),q_{\infty }}}=\Big(\sum_{v=-\infty }^{0}2^{-v\theta
q(0)}J(2^{v},u_{v})^{q(0)}\Big)^{\frac{1}{q(0)}}+\Big(\sum_{v=1}^{\infty
}2^{-v\theta q_{\infty }}J(2^{v},u_{v})^{q_{\infty }}\Big)^{\frac{1}{%
q_{\infty }}}<\infty .
\end{equation*}%
Moreover 
\begin{equation*}
\left\Vert f\right\Vert _{(A_{0},A_{1})_{\theta ,q(\cdot ),J}}\approx
\inf_{u_{v}}\left\Vert \left( J(2^{v},u_{v})\right) _{v}\right\Vert
_{\lambda ^{\theta ,q(0),q_{\infty }}},
\end{equation*}%
where the infimum is extended over all sequences $(u_{v})_{v}$ satisfying\ $%
\mathrm{\eqref{f-representation}}$.
\end{lemma}

\textbf{Proof.} Let $f\in (A_{0},A_{1})_{\theta ,q(\cdot ),J}$. Then we have
a representation%
\begin{equation*}
f=\int_{0}^{\infty }u(t)\frac{dt}{t},
\end{equation*}%
where $u(t)$ is measurable with values in $A_{0}\cap A_{1}$\ and%
\begin{equation*}
\left\Vert f\right\Vert _{(A_{0},A_{1})_{\theta ,q(\cdot ),J}}=\inf
\left\Vert t^{-\theta }J(t,u(t))\right\Vert _{L^{q(\cdot )}((0,\infty ),%
\frac{dt}{t})}<\infty .
\end{equation*}%
We set%
\begin{equation*}
u_{v}=\int_{2^{v}}^{2^{v+1}}u(t)\frac{dt}{t},\quad v\in \mathbb{Z}.
\end{equation*}%
Then we have%
\begin{equation*}
f=\sum_{v=-\infty }^{\infty }u_{v}.
\end{equation*}%
Let us prove that%
\begin{equation}
S(\{u_{v}\})=\Big(\sum_{v=-\infty }^{0}2^{-v\theta q(0)}\alpha _{v}^{q(0)}%
\Big)^{\frac{1}{q(0)}}+\Big(\sum_{v=1}^{\infty }2^{-v\theta q_{\infty
}}\alpha _{v}^{q_{\infty }}\Big)^{\frac{1}{q_{\infty }}}\lesssim \left\Vert
f\right\Vert _{(A_{0},A_{1})_{\theta ,q(\cdot ),J}},  \label{est-J}
\end{equation}%
with $\alpha _{v}=J(2^{v},u_{v})$, $v\in \mathbb{Z}$. We need only to prove
that%
\begin{equation*}
\sum_{v=1}^{\infty }2^{-v\theta q_{\infty }}\big(\frac{\alpha _{v}}{%
\left\Vert f\right\Vert _{(A_{0},A_{1})_{\theta ,q(\cdot ),J}}}\big)%
^{q_{\infty }}\lesssim 1,
\end{equation*}%
and%
\begin{equation*}
\sum_{v=-\infty }^{0}2^{-v\theta q(0)}\big(\frac{\alpha _{v}}{\left\Vert
f\right\Vert _{(A_{0},A_{1})_{\theta ,q(\cdot ),J}}}\big)^{q(0)}\lesssim 1.
\end{equation*}%
First let us prove that 
\begin{equation*}
2^{-v\theta q_{\infty }}\Big(\frac{\alpha _{v}}{\left\Vert f\right\Vert
_{(A_{0},A_{1})_{\theta ,q(\cdot ),J}}}\Big)^{q_{\infty }}\leq
\int_{2^{v}}^{2^{v+1}}\big(t^{-\theta }J(t,\frac{u(t)}{\left\Vert
f\right\Vert _{(A_{0},A_{1})_{\theta ,q(\cdot ),J}}})\big)^{q\left( t\right)
}\frac{dt}{t}+2^{-v}=\delta
\end{equation*}%
for any $v\in \mathbb{N}$. This claim can be reformulated as showing that%
\begin{equation*}
\Big(\delta ^{-\frac{1}{q_{\infty }}}2^{-v\theta }\frac{\alpha _{v}}{%
\left\Vert f\right\Vert _{(A_{0},A_{1})_{\theta ,q(\cdot ),J}}}\Big)%
^{q_{\infty }}\lesssim 1.
\end{equation*}%
Using the property $\mathrm{\eqref{Est-J}}$, we find that%
\begin{equation*}
\delta ^{-\frac{1}{q_{\infty }}}2^{-v\theta }\frac{\alpha _{v}}{\left\Vert
f\right\Vert _{(A_{0},A_{1})_{\theta ,q(\cdot ),J}}}\lesssim
\int_{2^{v}}^{2^{v+1}}\delta ^{-\frac{1}{q_{\infty }}}\tau ^{-\theta }J(\tau
,\frac{u(t)}{\left\Vert f\right\Vert _{(A_{0},A_{1})_{\theta ,q(\cdot ),J}}})%
\frac{d\tau }{\tau }.
\end{equation*}%
By Lemma \ref{DHHR-estimate},%
\begin{eqnarray*}
&&\Big(\int_{2^{v}}^{2^{v+1}}\delta ^{-\frac{1}{q_{\infty }}}\tau ^{-\theta
}J(\tau ,\frac{u(t)}{\left\Vert f\right\Vert _{(A_{0},A_{1})_{\theta
,q(\cdot ),J}}}\frac{d\tau }{\tau }\Big)^{q\left( t\right) } \\
&\lesssim &\int_{2^{v}}^{2^{v+1}}\delta ^{-\frac{q\left( \tau \right) }{%
q_{\infty }}}\big(\tau ^{-\theta }J(\tau ,\frac{u(t)}{\left\Vert
f\right\Vert _{(A_{0},A_{1})_{\theta ,q(\cdot ),J}}})\big)^{q\left( \tau
\right) }\frac{d\tau }{\tau }+1
\end{eqnarray*}%
for any $t\in \lbrack 2^{v},2^{v+1}]$. Since, $q$ is \emph{$\log $-}H\"{o}%
lder continuous at the infinity we find that%
\begin{equation}
\delta ^{-\frac{q\left( \tau \right) }{q_{\infty }}}\approx \delta ^{-1},%
\text{ \ \ }\tau \in \lbrack 2^{v},2^{v+1}],\quad v\in \mathbb{N}.
\label{Est-delta}
\end{equation}%
Therefore, from the definition of $\delta $, we find that the last integral
is dominated by a constant independent on $v\in \mathbb{N}$. Now, let us
prove that 
\begin{equation*}
2^{-v\theta q(0)}\Big(\frac{\alpha _{v}}{\left\Vert f\right\Vert
_{(A_{0},A_{1})_{\theta ,q(\cdot )}}}\Big)^{q(0)}\lesssim
\int_{2^{v}}^{2^{v+1}}\big(t^{-\theta }J(t,\frac{f}{\left\Vert f\right\Vert
_{(A_{0},A_{1})_{\theta ,q(\cdot ),J}}})\big)^{q\left( t\right) }\frac{dt}{t}%
+2^{v}=\delta
\end{equation*}%
for any $v\leq 0$. This claim can be reformulated as showing that%
\begin{equation*}
\Big(\delta ^{-\frac{1}{q(0)}}2^{-v\theta }\frac{\alpha _{v}}{\left\Vert
f\right\Vert _{(A_{0},A_{1})_{\theta ,q(\cdot ),J}}}\Big)^{q(0)}\lesssim 1.
\end{equation*}%
The property $\mathrm{\eqref{Est-J}}$, gives\ that%
\begin{equation*}
\delta ^{-\frac{1}{q(0)}}2^{-v\theta }\frac{\alpha _{v}}{\left\Vert
f\right\Vert _{(A_{0},A_{1})_{\theta ,q(\cdot ),J}}}\leq
\int_{2^{v}}^{2^{v+1}}\delta ^{-\frac{1}{q(0)}}\tau ^{-\theta }J(\tau ,\frac{%
u(t)}{\left\Vert f\right\Vert _{(A_{0},A_{1})_{\theta ,q(\cdot ),J}}})\frac{%
d\tau }{\tau }.
\end{equation*}%
Again by Lemma \ref{DHHR-estimate},%
\begin{eqnarray*}
&&\Big(\int_{2^{v}}^{2^{v+1}}\delta ^{-\frac{1}{q(0)}}\tau ^{-\theta }J(\tau
,\frac{u(t)}{\left\Vert f\right\Vert _{(A_{0},A_{1})_{\theta ,q(\cdot ),J}}})%
\frac{d\tau }{\tau }\Big)^{q\left( t\right) } \\
&\lesssim &\int_{2^{v}}^{2^{v+1}}\delta ^{-\frac{q\left( \tau \right) }{q(0)}%
}\big(\tau ^{-\theta }J(\tau ,\frac{u(t)}{\left\Vert f\right\Vert
_{(A_{0},A_{1})_{\theta ,q(\cdot ),J}}})\frac{d\tau }{\tau }\big)^{q\left(
\tau \right) }\frac{d\tau }{\tau }+1
\end{eqnarray*}%
for any $t\in \lbrack 2^{v},2^{v+1}]$\ and any $v\leq 0$. We use the
logarithmic decay\ condition at origin\ of $q$ to show that%
\begin{equation*}
\delta ^{-\frac{q\left( \tau \right) }{q(0)}}\approx \delta ^{-1},\text{ \ \ 
}\tau \in \lbrack 2^{v},2^{v+1}],\quad v\leq 0.
\end{equation*}%
Therefore and from the definition of $\delta $, we find that%
\begin{equation*}
\int_{2^{v}}^{2^{v+1}}\delta ^{-\frac{q\left( \tau \right) }{q(0)}}\big(\tau
^{-\theta }J(\tau ,\frac{u(t)}{\left\Vert f\right\Vert
_{(A_{0},A_{1})_{\theta ,q(\cdot ),J}}})\big)^{q\left( \tau \right) }\frac{%
d\tau }{\tau }\lesssim 1.
\end{equation*}%
Hence the left-hand side of $\mathrm{\eqref{est-J}}$ can be estimated by%
\begin{equation*}
c\int_{0}^{\infty }\big(t^{-\theta }J(t,\frac{u(t)}{\left\Vert f\right\Vert
_{(A_{0},A_{1})_{\theta ,q(\cdot ),J}}})\big)^{q\left( t\right) }\frac{dt}{t}%
+1.
\end{equation*}%
The first term is bounded since 
\begin{equation*}
\Big\|t^{-\theta }J(t,\frac{u(t)}{\left\Vert f\right\Vert
_{(A_{0},A_{1})_{\theta ,q(\cdot ),J}}})\Big\|_{L^{q(\cdot )}((0,\infty ),%
\frac{dt}{t})}\leq 1.
\end{equation*}%
Now in $\mathrm{\eqref{est-J}}$ taking the infimum, we conclude that%
\begin{equation*}
\inf_{u_{v}}\left\Vert \left( J(2^{v},u_{v})\right) _{v}\right\Vert
_{\lambda ^{\theta ,q(0),q_{\infty }}}\lesssim \left\Vert f\right\Vert
_{(A_{0},A_{1})_{\theta ,q(\cdot ),J}}.
\end{equation*}%
Conversely, assume that 
\begin{equation*}
f=\sum_{v=-\infty }^{\infty }u_{v},
\end{equation*}%
and 
\begin{equation*}
S(\{u_{v}\})<\infty .
\end{equation*}%
Let us prove that%
\begin{equation*}
\left\Vert f\right\Vert _{(A_{0},A_{1})_{\theta ,q(\cdot ),J}}\lesssim
\inf_{w_{v}}S(w_{v})=S,
\end{equation*}%
where the infimum is taking over all sequences $\{w_{v}\}$ satisfying $%
\mathrm{\eqref{f-representation}}$. Choose%
\begin{equation*}
u(t)=\frac{u_{v}}{\log 2},\quad t\in \lbrack 2^{v},2^{v+1}].
\end{equation*}%
Then $f=\int_{0}^{\infty }u(t)\frac{dt}{t}$. This claim can be reformulated
as showing that%
\begin{equation*}
\int_{0}^{\infty }\big(t^{-\theta }J(t,\frac{u(t)}{S})\big)^{q\left(
t\right) }\frac{dt}{t}\lesssim 1.
\end{equation*}%
Now our estimate clearly follows from the inequalities%
\begin{equation*}
\int_{2^{v}}^{2^{v+1}}\big(t^{-\theta }J(t,\frac{u(t)}{S})\big)^{q\left(
t\right) }\frac{dt}{t}\lesssim 2^{-v\theta q_{\infty }}\big(\frac{\alpha _{v}%
}{S}\big)^{q_{\infty }}+2^{-v}=\delta
\end{equation*}%
for any $v\in \mathbb{N}$ and 
\begin{equation*}
\int_{2^{v}}^{2^{v+1}}\big(t^{-\theta }J(t,\frac{u(t)}{S})\big)^{q\left(
t\right) }\frac{dt}{t}\lesssim 2^{-v\theta q(0)}\big(\frac{\alpha _{v}}{S}%
\big)^{q(0)}+2^{v}=\delta
\end{equation*}%
for any $v\leq 0$. The first claim can be reformulated as showing\ that%
\begin{equation*}
\int_{2^{v}}^{2^{v+1}}\big(\delta ^{-\frac{1}{q(t)}}t^{-\theta }J(t,\frac{%
u(t)}{S})\big)^{q\left( t\right) }\frac{dt}{t}\lesssim 1.
\end{equation*}%
We need only to show that%
\begin{equation*}
\delta ^{-\frac{1}{q(t)}}t^{-\theta }J(t,\frac{u(t)}{S})\lesssim 1
\end{equation*}%
for any $v\in \mathbb{N}$ and any $t\in \lbrack 2^{v},2^{1+v}]$. The
left-hand side is bounded by%
\begin{equation*}
\delta ^{-\frac{1}{q(t)}}2^{-\theta v}J(2^{v},\frac{u(t)}{S}).
\end{equation*}%
From $\mathrm{\eqref{Est-delta}}$ we find that 
\begin{equation*}
\delta ^{-\frac{1}{q(t)}}2^{-\theta v}J(2^{v},\frac{u(t)}{S}%
,A_{0},A_{1})\lesssim \delta ^{-\frac{1}{q_{\infty }}}2^{-\theta v}J(2^{v},%
\frac{u_{v}}{S})\leq 1
\end{equation*}%
for any $v\in \mathbb{N}$. Similarly we estimate the second claim.

We conclude that%
\begin{equation*}
\inf_{u_{v}}\left\Vert \left( J(2^{v},u_{v})\right) _{v}\right\Vert
_{\lambda ^{\theta ,q(0),q_{\infty }}}\lesssim \left\Vert f\right\Vert
_{(A_{0},A_{1})_{\theta ,q(\cdot ),J}}.
\end{equation*}%
\quad $\square $

We shall prove that the spaces generated by the $K$-and $J$-methods are the
same.

\begin{theorem}
\label{Equi-J-and-K}Let $(A_{0},A_{1})$\ be a compatible couple.\ Let $%
\theta \in (0,1)$\textit{\ and }$q\in \mathcal{P}(\mathbb{R})$ \textit{be
log-H\"{o}lder continuous both at the origin and at infinity. Then }%
\begin{equation*}
(A_{0},A_{1})_{\theta ,q(\cdot ),J}=(A_{0},A_{1})_{\theta ,q(\cdot )},
\end{equation*}%
with equivalence of norms.
\end{theorem}

\textbf{Proof.} Let $f\in (A_{0},A_{1})_{\theta ,q(\cdot ),J}$ with $%
f=\int_{0}^{\infty }u(s)\frac{ds}{s}$, where $u(t)$ is measurable with
values in $A_{0}\cap A_{1}$. By $\mathrm{\eqref{Est-J1}}$ we have 
\begin{equation*}
K(t,f)\leq \int_{0}^{\infty }K(t,u(s))\frac{ds}{s}\leq \int_{0}^{\infty
}\min (1,\frac{t}{s})J(s,u(s))\frac{ds}{s}.
\end{equation*}%
Applying Hardy inequality, Lemma \ref{lq-inequality}, we get%
\begin{equation*}
\left\Vert f\right\Vert _{(A_{0},A_{1})_{\theta ,q(\cdot )}}\lesssim
\left\Vert f\right\Vert _{(A_{0},A_{1})_{\theta ,q(\cdot ),J}}.
\end{equation*}%
For the converse inequality, Lemma 3.3.2 of \cite{BL76}, and using Theorem %
\ref{embedding}, implies the existence of a representation 
\begin{equation*}
f=\sum_{v=-\infty }^{\infty }u_{v},
\end{equation*}%
such that 
\begin{equation*}
J(2^{v},u_{v})\leq (\gamma +\varepsilon )K(2^{v},f)
\end{equation*}%
for any $v\in \mathbb{Z}$, $\varepsilon >0$ and $\gamma $ is a universal
constant less than or equal $3$. By Lemmas \ref{K-property2} and \ref%
{J-property2} we get%
\begin{equation*}
\left\Vert f\right\Vert _{(A_{0},A_{1})_{\theta ,q(\cdot ),J}}\lesssim
\left\Vert f\right\Vert _{(A_{0},A_{1})_{\theta ,q(\cdot )}}.
\end{equation*}%
This completes the proof of this theorem.$\quad \square $

\begin{theorem}
\label{Density}Let $(A_{0},A_{1})$\ be a compatible couple.\ Let $\theta \in
(0,1)$\textit{\ and }$q\in \mathcal{P}(\mathbb{R})$ \textit{be log-H\"{o}%
lder continuous both at the origin and at infinity. Then }$A_{0}\cap A_{1}$
is dense in $(A_{0},A_{1})_{\theta ,q(\cdot )}$.
\end{theorem}

\textbf{Proof.} Let $f\in (A_{0},A_{1})_{\theta ,q(\cdot )}$. From Theorem %
\ref{Equi-J-and-K} we have%
\begin{equation*}
f=\sum_{v=-\infty }^{\infty }u_{v},
\end{equation*}%
where $u_{v}$, $v\in \mathbb{Z}$ is measurable with values in $A_{0}\cap
A_{1}$ and%
\begin{equation*}
\left\Vert \left( J(2^{v},u_{v})\right) _{v}\right\Vert _{\lambda ^{\theta
,q(0),q_{\infty }}}<\infty .
\end{equation*}%
Then%
\begin{eqnarray*}
&&\Big\|f-\sum_{|v|\leq N}u_{v}\Big\|_{(A_{0},A_{1})_{\theta ,q(\cdot )}} \\
&\leq &\Big(\sum_{v=N}^{\infty }2^{-v\theta q_{\infty
}}J(2^{v},u_{v})^{q_{\infty }}\Big)^{\frac{1}{q_{\infty }}}+\Big(%
\sum_{v=-\infty }^{-N}2^{-v\theta q(0)}J(2^{v},u_{v})^{q(0)}\Big)^{\frac{1}{%
q(0)}}.
\end{eqnarray*}%
Therefore,%
\begin{equation*}
\Big\|f-\sum_{|v|\leq N}u_{v}\Big\|_{(A_{0},A_{1})_{\theta ,q(\cdot )}},
\end{equation*}%
which tends to zero if $N\longrightarrow \infty $.$\quad \square $

\begin{definition}
Let $\theta \in \lbrack 0,1]$. Let $(A_{0},A_{1})$ be a compatible couple of
normed vector spaces. Suppose that $X$ is an intermediate space with respect
to $(A_{0},A_{1})$. Then we say that \newline
(i) $X$ is of class $\mathscr{C}_{K}(\theta ;A_{0},A_{1})$ if $%
K(t,f;A_{0},A_{1})\leq Ct^{\theta }\left\Vert f\right\Vert _{X},\quad f\in
X; $\newline
(ii) $X$ is of class $\mathscr{C}_{J}(\theta ;A_{0},A_{1})$ if $\left\Vert
f\right\Vert _{X}\leq Ct^{-\theta }J(t,f;A_{0},A_{1}),\quad f\in A_{0}\cap
A_{1}.$\newline
(iii) We say that $X$ is of class $\mathscr{C}(\theta ;A_{0},A_{1})$ if $X$
is of class $\mathscr{C}_{K}(\theta ;A_{0},A_{1})$ and of class\ $\mathscr{C}%
_{J}(\theta ;A_{0},A_{1})$.
\end{definition}

Let $q\in \mathcal{P}(\mathbb{R})$\textit{. }From \cite[Theorem 3.5.2]{BL76}
and Proposition 3 we see that $(A_{0},A_{1})_{\theta ,q(\cdot )}$ is of
class $\mathscr{C}(\theta ;A_{0},A_{1})$ if $\theta \in (0,1)$.

We are now ready to prove the reiteration theorem, which is one of the most
important general results in interpolation theory.

\begin{theorem}
\label{The reiteration theorem}Let $q\in \mathcal{P}(\mathbb{R})$ \textit{be
log-H\"{o}lder continuous both at the origin and at infinity. }Let $%
(A_{0},A_{1})$ and $(X_{0},X_{1})$ be two compatible couples of normed
linear spaces, and assume that $X_{i}$ ($i=0,1$) are complete and of class $%
\mathscr{C}(\theta _{i};A_{0},A_{1})$, where $\theta _{0},\theta _{1}\in
\lbrack 0,1]$ and $\theta _{0}\neq \theta _{1}$. Put%
\begin{equation*}
\theta =(1-\eta )\theta _{0}+\eta \theta _{1},\quad \eta \in (0,1).
\end{equation*}%
Then%
\begin{equation*}
(A_{0},A_{1})_{\theta ,q(\cdot )}=(X_{0},X_{1})_{\eta ,q(\cdot )}
\end{equation*}%
with equivalence of norms. In particular, if $\theta _{0},\theta _{1}\in
(0,1)$, $q_{0},q_{1}\in \mathcal{P}(\mathbb{R})$ \textit{are log-H\"{o}lder
continuous both at the origin and at infinity} and $(A_{0},A_{1})_{\theta
_{i},q_{i}(\cdot )}$ are complete then 
\begin{equation*}
\left( (A_{0},A_{1})_{\theta _{0},q_{0}(\cdot )},(A_{0},A_{1})_{\theta
_{1},q_{1}(\cdot )}\right) _{\eta ,q(\cdot )}=(A_{0},A_{1})_{\theta ,q(\cdot
)}
\end{equation*}%
where%
\begin{equation*}
\frac{1}{q(\cdot )}=\frac{\theta _{0}}{q_{0}(\cdot )}+\frac{\theta _{1}}{%
q_{1}(\cdot )}.
\end{equation*}
\end{theorem}

\textbf{Proof. }We will do the proof in two steps.

\textit{Step 1.} Let us prove that%
\begin{equation}
(X_{0},X_{1})_{\eta ,q(\cdot )}\hookrightarrow (A_{0},A_{1})_{\theta
,q(\cdot )}.  \label{First-emb}
\end{equation}%
Let $f\in (X_{0},X_{1})_{\eta ,q(\cdot )}$. Then%
\begin{equation*}
f=f_{0}+f_{1},\quad f_{0}\in X_{0},f_{1}\in X_{1}.
\end{equation*}%
Since $X_{i}$ ($i=0,1$) are of class $\mathscr{C}(\theta _{i};A_{0},A_{1})$
we have%
\begin{eqnarray*}
K(t,f;A_{0},A_{1}) &\leq &K(t,f_{0};A_{0},A_{1})+K(t,f_{1};A_{0},A_{1})\leq
c(t^{\theta _{0}}\left\Vert f_{0}\right\Vert _{X_{0}}+t^{\theta
_{1}}\left\Vert f_{1}\right\Vert _{X_{1}}) \\
&\leq &ct^{\theta _{0}}K(t^{\theta _{1}-\theta _{0}},f;X_{0},X_{1}).
\end{eqnarray*}%
Therefore, from Theorem \ref{K-property2}, we get 
\begin{eqnarray*}
\left\Vert f\right\Vert _{(A_{0},A_{1})_{\theta ,q(\cdot )}} &\lesssim &\Big(%
\int_{0}^{1}t^{(\theta _{0}-\theta )q(0)}K(t^{\theta _{1}-\theta
_{0}},f;X_{0},X_{1})^{q(0)}\frac{dt}{t}\Big)^{\frac{1}{q(0)}} \\
&&+\Big(\int_{1}^{\infty }t^{(\theta _{0}-\theta )q_{\infty }}K(t^{\theta
_{1}-\theta _{0}},f;X_{0},X_{1})^{q_{\infty }}\frac{dt}{t}\Big)^{\frac{1}{%
q_{\infty }}}.
\end{eqnarray*}%
Putting\ $s=t^{\theta _{1}-\theta _{0}}$ and observing that $\eta =\frac{%
\theta -\theta _{0}}{\theta _{1}-\theta _{0}}$ we find that%
\begin{equation*}
\left\Vert f\right\Vert _{(A_{0},A_{1})_{\theta ,q(\cdot )}}\lesssim
\left\Vert f\right\Vert _{(X_{0},X_{1})_{\eta ,q(\cdot )}},
\end{equation*}%
which gives $\mathrm{\eqref{First-emb}}$.

\textit{Step 2.} Let us prove that%
\begin{equation}
(A_{0},A_{1})_{\theta ,q(\cdot )}\hookrightarrow (X_{0},X_{1})_{\eta
,q(\cdot )}.  \label{Second-emb}
\end{equation}%
Assume that $f\in (A_{0},A_{1})_{\theta ,q(\cdot )}$. We choose a
representation 
\begin{equation*}
f=\sum_{v=-\infty }^{\infty }u_{v},
\end{equation*}%
where $u_{v}$, $v\in \mathbb{Z}$ is measurable with values in $A_{0}\cap
A_{1}$ and%
\begin{equation*}
\left\Vert \left( J(2^{v},u_{v})\right) _{v}\right\Vert _{\lambda ^{\theta
,q(0),q_{\infty }}}<\infty .
\end{equation*}%
Applying $\mathrm{\eqref{Est-J1}}$, and that $X_{i}$ ($i=0,1$) are of class $%
\mathscr{C}(\theta _{i};A_{0},A_{1})$ we get for any $j\in \mathbb{Z}$,%
\begin{eqnarray*}
&&2^{(\theta _{0}-\theta )j}K(2^{(\theta _{1}-\theta _{0})j},f;X_{0},X_{1})
\\
&\leq &2^{(\theta _{0}-\theta )j}\sum_{v=-\infty }^{\infty }K(2^{(\theta
_{1}-\theta _{0})j},u_{v};X_{0},X_{1}) \\
&\leq &2^{(\theta _{0}-\theta )j}\sum_{v=-\infty }^{\infty }\min \left(
1,2^{(j-v)(\theta _{1}-\theta _{0})}\right) J(2^{v(\theta _{1}-\theta
_{0})},u_{v};X_{0},X_{1}) \\
&\leq &2^{-\theta j}\sum_{v=-\infty }^{\infty }\min \left( 2^{(j-v)\theta
_{0}},2^{(j-v)\theta _{1}}\right) J(2^{v},u_{v};A_{0},A_{1}).
\end{eqnarray*}%
The last term can be rewritten us%
\begin{equation}
\sum_{v=-\infty }^{j}2^{(j-v)(\theta _{0}-\theta )}2^{-v\theta
}J(2^{v},u_{v};A_{0},A_{1})+\sum_{v=j+1}^{\infty }2^{(j-v)(\theta
_{1}-\theta )}2^{-v\theta }J(2^{v},u_{v};A_{0},A_{1})  \label{sum-equi}
\end{equation}%
for any $j\in \mathbb{Z}$. We treat the case where\ $j\geq 0$. The first sum
can be rewritten us%
\begin{eqnarray*}
&&\sum_{v=-\infty }^{0}2^{(j-v)(\theta _{0}-\theta )}2^{-v\theta
}J(2^{v},u_{v};A_{0},A_{1})+\sum_{v=1}^{j}2^{(j-v)(\theta _{0}-\theta
)}2^{-v\theta }J(2^{v},u_{v};A_{0},A_{1}) \\
&\lesssim &2^{j(\theta _{0}-\theta )}\sup_{v\leq 0}(2^{-v\theta
}J(2^{v},u_{v};A_{0},A_{1}))+\sum_{v=1}^{j}2^{(j-v)(\theta _{0}-\theta
)}2^{-v\theta }J(2^{v},u_{v};A_{0},A_{1}).
\end{eqnarray*}%
Applying Lemma \ref{Hardy-inequality} we get%
\begin{equation*}
\left\Vert \left( 2^{(\theta _{0}-\theta )j}K(2^{(\theta _{1}-\theta
_{0})j},f;X_{0},X_{1})\right) _{j\geq 1}\right\Vert _{\lambda ^{\theta
,q(0),q_{\infty }}}\leq \left\Vert \left( J(2^{v},u_{v})\right)
_{v}\right\Vert _{\lambda ^{\theta ,q(0),q_{\infty }}}.
\end{equation*}%
Now if $j\leq 0$, the second sum of $\mathrm{\eqref{sum-equi}}$ can be
rewritten us%
\begin{eqnarray*}
&&\sum_{v=j+1}^{0}2^{(j-v)(\theta _{1}-\theta )}2^{-v\theta
}J(2^{v},u_{v};A_{0},A_{1})+\sum_{v=1}^{\infty }2^{(j-v)(\theta _{1}-\theta
)}2^{-v\theta }J(2^{v},u_{v};A_{0},A_{1}) \\
&\leq &\sum_{v=j+1}^{0}2^{(j-v)(\theta _{1}-\theta )}2^{-v\theta
}J(2^{v},u_{v};A_{0},A_{1})+2^{j(\theta _{1}-\theta )}\sup_{v\geq 1}\left(
2^{-v\theta }J(2^{v},u_{v};A_{0},A_{1})\right) .
\end{eqnarray*}%
Applying again Lemma \ref{Hardy-inequality} we get%
\begin{equation*}
\left\Vert \left( 2^{(\theta _{0}-\theta )j}K(2^{(\theta _{1}-\theta
_{0})j},f;X_{0},X_{1})\right) _{j\leq 0}\right\Vert _{\lambda ^{\theta
,q(0),q_{\infty }}}\leq \left\Vert \left( J(2^{v},u_{v})\right)
_{v}\right\Vert _{\lambda ^{\theta ,q(0),q_{\infty }}}.
\end{equation*}%
This prove the embedding $\mathrm{\eqref{Second-emb}}$ by taking the infimum
in view of the Theorem \ref{J-property2} and the fact\ that%
\begin{equation*}
\left\Vert f\right\Vert _{(X_{0},X_{1})_{\eta ,q(\cdot )}}\approx \Big(%
\sum_{j=-\infty }^{0}2^{(\theta _{0}-\theta )q(0)}\alpha _{j}^{q(0)}\Big)^{%
\frac{1}{q(0)}}+\Big(\sum_{j=1}^{\infty }2^{(\theta _{0}-\theta )q_{\infty
}}\alpha _{j}^{q_{\infty }}\Big)^{\frac{1}{q_{\infty }}},
\end{equation*}%
where%
\begin{equation*}
\alpha _{j}=K(2^{(\theta _{1}-\theta _{0})j},f;X_{0},X_{1}),\quad j\in 
\mathbb{Z}.
\end{equation*}%
This completes the proof of Theorem \ref{The reiteration theorem}.$\quad
\square $

\section{Application}

In this section, we give a simple application of the results of the previous
sections. We will present various real interpolation formulas in Besov
spaces with variable indices. The symbol $\mathcal{S}(\mathbb{R}^{n})$ is
used in place of the set of all Schwartz functions on $\mathbb{R}^{n}$. We
denote by $\mathcal{S}^{\prime }(\mathbb{R}^{n})$ the dual space of all
tempered distributions on $\mathbb{R}^{n}$. The Fourier transform of a
Schwartz function $f$ is denoted by $\mathcal{F}f$. To define the variable
Besov spaces, we first need the concept of a smooth dyadic resolution of
unity. Let $\Psi $\ be a function\ in $\mathcal{S}(%
\mathbb{R}
^{n})$\ satisfying $\Psi (x)=1$\ for\ $\left\vert x\right\vert \leq 1$\ and\ 
$\Psi (x)=0$\ for\ $\left\vert x\right\vert \geq 2$.\ We define $\varphi
_{0} $ and $\varphi _{1}$ by $\mathcal{F}\varphi _{0}(x)=\Psi (x)$, $%
\mathcal{F}\varphi _{1}(x)=\Psi (x)-\Psi (2x)$\ and 
\begin{equation*}
\mathcal{F}\varphi _{j}(x)=\mathcal{F}\varphi _{1}(2^{-j}x)\quad \text{%
\textit{for}}\quad j=2,3,....
\end{equation*}%
Then $\{\mathcal{F}\varphi _{j}\}_{j\in \mathbb{N}_{0}}$\ is a smooth dyadic
resolution of unity, $\sum_{j=0}^{\infty }\mathcal{F}\varphi _{j}(x)=1$ for
all $x\in \mathbb{R}^{n}$.\ Thus we obtain the Littlewood-Paley
decomposition $f=\sum_{j=0}^{\infty }\varphi _{j}\ast f$ of all $f\in 
\mathcal{S}^{\prime }(%
\mathbb{R}
^{n})$ $($convergence in $\mathcal{S}^{\prime }(%
\mathbb{R}
^{n}))$.

Let $p,q\in \mathcal{P}(\mathbb{R}^{n})$. The mixed Lebesgue-sequence space $%
\ell _{>}^{q(\cdot )}(L^{p(\cdot )})$ is defined on sequences of $L^{p(\cdot
)}$-functions by the modular%
\begin{equation*}
\varrho _{\ell _{>}^{q(\cdot )}(L^{p\left( \cdot \right)
})}((f_{v})_{v}):=\sum\limits_{v=1}^{\infty }\inf \Big\{\lambda
_{v}>0:\varrho _{p(\cdot )}\Big(\frac{f_{v}}{\lambda _{v}^{1/q(\cdot )}}\Big)%
\leq 1\Big\}.
\end{equation*}%
The (quasi)-norm is defined from this as usual:%
\begin{equation}
\left\Vert \left( f_{v}\right) _{v}\right\Vert _{\ell _{>}^{q(\cdot
)}(L^{p\left( \cdot \right) })}:=\inf \Big\{\mu >0:\varrho _{\ell ^{q(\cdot
)}(L^{p(\cdot )})}\Big(\frac{1}{\mu }(f_{v})_{v}\Big)\leq 1\Big\}.
\label{mixed-norm}
\end{equation}%
If $q^{+}<\infty $, then we can replace $\mathrm{\eqref{mixed-norm}}$ by the
simpler expression $\varrho _{\ell _{>}^{q(\cdot )}(L^{p(\cdot
)})}((f_{v})_{v}):=\sum\limits_{v=1}^{\infty }\left\Vert |f_{v}|^{q(\cdot
)}\right\Vert _{\frac{p(\cdot )}{q(\cdot )}}$. The case $p:=\infty $ can be
included by replacing the last modular by $\varrho _{\ell _{>}^{q(\cdot
)}(L^{\infty })}((f_{v})_{v}):=\sum\limits_{v=1}^{\infty }\big\|\left\vert
f_{v}\right\vert ^{q(\cdot )}\big\|_{\infty }$.

We define the following class of variable exponents $\mathcal{P}^{\mathrm{log%
}}(\mathbb{R}^{n}):=\big\{p\in \mathcal{P}:\frac{1}{p}\in C^{\log }\big\}$,
were introduced in $\mathrm{\cite[Section \ 2]{DHHMS}}$. We define $%
1/p_{\infty }:=\lim_{|x|\rightarrow \infty }1/p(x)$\ and we use the
convention $\frac{1}{\infty }=0$. Note that although $\frac{1}{p}$ is
bounded, the variable exponent $p$ itself can be unbounded.

We state the definition of the spaces $B_{p(\cdot ),q(\cdot )}^{s(\cdot )}$,
which introduced and investigated in \cite{AH}.

\begin{definition}
\label{11}\textit{Let }$\left\{ \mathcal{F}\varphi _{j}\right\}
_{j=0}^{\infty }$\textit{\ be a resolution of unity}, $s:\mathbb{R}%
^{n}\rightarrow \mathbb{R}$ and $p,q\in \mathcal{P}(\mathbb{R}^{n})$.\textit{%
\ The Besov space }$B_{p(\cdot ),q(\cdot )}^{s(\cdot )}$\textit{\ consists
of all distributions }$f\in \mathcal{S}^{\prime }(%
\mathbb{R}
^{n})$\textit{\ such that}%
\begin{equation*}
\left\Vert f\right\Vert _{B_{p(\cdot ),q(\cdot )}^{s(\cdot )}}:=\left\Vert
(2^{js(\cdot )}\varphi _{j}\ast f)_{j}\right\Vert _{\ell _{>}^{q(\cdot
)}(L^{p\left( \cdot \right) })}<\infty .
\end{equation*}
\end{definition}

Taking $s\in \mathbb{R}$ and $q\in (0,\infty ]$ as constants we derive the
spaces $B_{p(\cdot ),q}^{s}$ studied by Xu in \cite{Xu08}. We refer the
reader to the recent papers \cite{D3}, \cite{KV121},\cite{KV122} and\ \cite%
{IN14} for further details, historical remarks and more references on these
function spaces. For any $p,q\in \mathcal{P}_{0}^{\log }(\mathbb{R}^{n})$
and $s\in C_{\text{loc}}^{\log }$, the space $B_{p(\cdot ),q(\cdot
)}^{s(\cdot )}$ does not depend on the chosen smooth dyadic resolution of%
\textit{\ }unity $\{\mathcal{F}\varphi _{j}\}_{j\in \mathbb{N}_{0}}$ (in the
sense of\ equivalent quasi-norms) and 
\begin{equation*}
\mathcal{S}(\mathbb{R}^{n})\hookrightarrow B_{p(\cdot ),q(\cdot )}^{s(\cdot
)}\hookrightarrow \mathcal{S}^{\prime }(\mathbb{R}^{n}).
\end{equation*}%
Moreover, if $p,q,s$ are constants, we re-obtain the usual Besov spaces $%
B_{p,q}^{s}$, studied in detail in \cite{Pe76}, \cite{SiRu96}, \cite{T1}, 
\cite{T2} and \cite{T31}.

Applying Lemma \ref{K-property2} and using the same arguments of \cite[%
Theorem 3.1]{AH14} we obtain.

\begin{theorem}
\label{real-int-Besov}Let $\theta \in (0,1)$\textit{\ and }$q\in \mathcal{P}(%
\mathbb{R})$ \textit{be log-H\"{o}lder continuous both at the origin and at
infinity\ with }$q(0)=q_{\infty }$. Let $p,q_{0},q_{1}\in \mathcal{P}^{\log
}(\mathbb{R}^{n})$ and $\alpha _{0},\alpha _{1}\in C_{\mathrm{loc}}^{\log }$%
. If $0\neq \alpha _{0}-\alpha _{1}$ is constant, then%
\begin{equation*}
(B_{p(\cdot ),q_{0}(\cdot )}^{\alpha _{0}(\cdot )},B_{p(\cdot ),q_{1}(\cdot
)}^{\alpha _{1}(\cdot )})_{\theta ,q(\cdot )}=B_{p(\cdot ),q(0)}^{\alpha
(\cdot )}
\end{equation*}%
with $\alpha (\cdot )=(1-\theta )\alpha _{0}(\cdot )+\theta \alpha
_{1}(\cdot )$. Moerover%
\begin{equation*}
(B_{p(\cdot ),r_{0}}^{\alpha _{0}(\cdot )},B_{p(\cdot ),r_{1}}^{\alpha
_{1}(\cdot )})_{\theta ,q(\cdot )}=B_{p(\cdot ),q(0)}^{\alpha (\cdot )},
\end{equation*}%
with $r_{0},r_{1}\in \lbrack 1,\infty ]$\ and 
\begin{equation*}
\frac{1}{q(0)}=\frac{1-\theta }{r_{0}}+\frac{\theta }{r_{1}}.
\end{equation*}
\end{theorem}

Now we present some interpolation results in variable exponent Lorentz
spaces $\mathcal{L}^{p(\cdot ),q(\cdot )}(\mathbb{R}^{n})$ introduced by 
\cite{EKS08}.

\begin{definition}
If $f$ is a measurable function on $\mathbb{R}^{n}$, we define the
non-increasing rearrangement of $f$ through%
\begin{equation*}
f^{\ast }(t)=\sup \{\lambda >0:m_{f}(\lambda )>t\}
\end{equation*}%
where $m_{f}$ is the distribution function of $f$.
\end{definition}

\begin{definition}
Let $p,q\in \mathcal{P}(\mathbb{R})$. By $\mathcal{L}^{p(\cdot ),q(\cdot )}(%
\mathbb{R}^{n})$\ we denote the space of functions f on $\mathbb{R}^{n}$
such that 
\begin{equation*}
\left\Vert f\right\Vert _{\mathcal{L}^{p(\cdot ),q(\cdot )}(\mathbb{R}^{n})}=%
\big\|t^{\frac{1}{p(t)}-\frac{1}{q(t)}}f^{\ast }(t)\big\|_{L^{q(\cdot )}(%
\mathbb{[}0,\infty ))}<\infty .
\end{equation*}
\end{definition}

We refer to the recent\ paper \cite{EKS08} for further details on these
scales of spaces. We present an equivalent quasi-norm for the space $%
\mathcal{L}^{p(\cdot ),q(\cdot )}(\mathbb{R}^{n})$, where the proof is quite
similar to that for Lemma \ref{K-property2}.

\begin{lemma}
\label{Lorentz-property}\textit{Let }$p,q\in \mathcal{P}(\mathbb{R})$ 
\textit{be log-H\"{o}lder continuous both at the origin and at infinity. Then%
}%
\begin{equation*}
\left\Vert f\right\Vert _{\mathcal{L}^{p(\cdot ),q(\cdot )}(\mathbb{R}%
^{n})}\approx \Big(\sum_{v=-\infty }^{0}2^{-v\frac{q(0)}{p(0)}}\left(
f^{\ast }(2^{v})\right) ^{q(0)}\Big)^{\frac{1}{q(0)}}+\Big(%
\sum_{v=1}^{\infty }2^{-v\frac{q_{\infty }}{p_{\infty }}}\left( f^{\ast
}(2^{v})\right) ^{q_{\infty }}\Big)^{\frac{1}{q_{\infty }}}.
\end{equation*}%
Moreover,%
\begin{equation*}
\left\Vert f\right\Vert _{\mathcal{L}^{p(\cdot ),q(\cdot )}(\mathbb{R}%
^{n})}\approx \Big(\int_{0}^{1}t^{-\frac{q(0)}{p(0)}}\left( f^{\ast
}(t)\right) ^{q(0)}\frac{dt}{t}\Big)^{\frac{1}{q(0)}}+\Big(\int_{1}^{\infty
}t^{-\frac{q_{\infty }}{p_{\infty }}}\left( f^{\ast }(t)\right) ^{q_{\infty
}}\frac{dt}{t}\Big)^{\frac{1}{q_{\infty }}}.
\end{equation*}
\end{lemma}

Applying this lemma and \cite[Theorem 5.2.1]{BL76} we obtain.

\begin{theorem}
\label{real-int-Lorentz}Let $\theta \in (0,1)$\textit{\ and }$q\in \mathcal{P%
}(\mathbb{R})$ \textit{be log-H\"{o}lder continuous both at the origin and
at infinity\ with }$q(0)=q_{\infty }$. Then%
\begin{equation*}
(L^{1},L_{\infty })_{\theta ,q(\cdot )}=\mathcal{L}^{p,q(\cdot )}(\mathbb{R}%
^{n}),
\end{equation*}%
with $p=\frac{1}{1-\theta }$.
\end{theorem}

Douadi Drihem

M'sila University, Department of Mathematics,

Laboratory of Functional Analysis and Geometry of Spaces,

P.O. Box 166, M'sila 28000, Algeria,

e-mail: \texttt{\ douadidr@yahoo.fr}

\end{document}